\documentclass[12pt,twoside,a4paper]{amsart}
\UseRawInputEncoding
\usepackage{amsmath,amssymb}
\usepackage{cite}
\usepackage{graphicx}
\usepackage{txfonts}
\usepackage{mathrsfs}
\newcommand{\ds}{\displaystyle}
\newcommand{\K}{J.Kaczorowski and K.Wiertelak}

\allowdisplaybreaks[4]

\theoremstyle{definition}

\theoremstyle{theorem}

\begin{document}

\begin{center}
\textbf{\large On the solution of the Volterra integral equation of second type for the error term in an asymptotic formula for arithmetical functions}
\\
\vspace{0.9cm}
by
\\
\vspace{0.35cm}
Hideto  IWATA
\end{center}

\vspace{0.5cm}

\textbf{Abstract.}\hspace{0.1cm}In 2010, {\K} considered the Volterra integral equation of second type for the remainder term in the asymptotic formula for the Euler totient function. The author found that the consideration made by them holds for other remainder terms in the asymptotic formula having  certain common properties. In this paper, we first consider the pair of complex-valued arithmetical functions $(a(n),b(n))$ satisfying  $b(n) \hspace{-0.1cm} = \sum_{d|n} a(d)n/d$. We prove that the solution of the Volterra integral equation of second type for the error term in the asymptotic formula for $b(n)$ can be obtained when $a(n)$ satisfies some special condition.

\vspace{1cm}
\begin{center}
\textbf{1.Introduction and the statement of the main result}
\end{center}
\vspace{0.1cm}

Let $\varphi(n)$ denote the Euler totient function and let
\vspace{-0.01cm}
\begin{equation}
   E(x) = \sum_{n \leq x} \varphi(n) - \frac{3}{\pi^2}x^2   \tag{1.1}
\end{equation}be the associated error term. {\K} studied the following Volterra integral equation of second type for $E(x)$ (see~\cite{Kac and Wie}) :
\vspace{-0.01cm}
\begin{equation}
   F(x) - \int_{0}^\infty K(x,t)F(t)dt = E(x)  \quad (x \geq 1), \tag{1.2}
\end{equation}where $F(x)$ is the unknown function and the kernel $K(x,t)$ is defined as follows:
\footnote[0]
{\textit{2010 Mathematics Subject Classification. Primary}
 45D05;11A25;11N37

\textit{Key words and phrases}
: Volterra integral equation of second type,
the remainder term in the asymptotic formula,
arithmetical function.
}
\begin{equation*}
   K(x,t) = \begin{cases}
                 1/t \quad (0 < t \leq x),
                 \\
                  0  \quad (1 \leq x < t). 
              \end{cases}
\end{equation*}

The equation (1.2) can be solved explicitly. Let us put
\begin{equation*}
   f(x) = -\sum_{n=1}^\infty \frac{\mu(n)}{n}\left\{ \frac{x}{n} \right\}   \tag{1.3}
\end{equation*}for every $x \geq 0$, where $\mu(n)$ denotes the M\"{o}bius function and $\{x\}=x-[x]$ is the fractional part of a real number $x$. Then the general solution of (1.2) is
\vspace{-0.01cm}
\begin{equation*}
   F(x) = (f(x)+A)x,   \tag{1.4}
\end{equation*}where $A$ is an arbitrary constant. (In the paper~\cite{Kac and Wie}, $F(x)=xf(x)$ is claimed to be the unique solution of the integral equation (1.2), but this uniqueness does not hold even assuming the initial value condition at $x=0$ . Probably, the term $Ax$ is missing to give the general solution.)

\vspace{0.1cm}

In the present article we generalize the observation of {\K} by slightly changing the settings, and propose the following conjecture.

\textsc{\large Conjecture.}  \textit{Let $\{ a(n) \}$ be a complex-valued arithmetical function satisfying some suitable conditions, and $\{ b(n) \}$ be the arithmetical function given by}
                                      \[ \ds b(n) = \sum_{d|n} a(d)\frac{n}{d}.\]\textit{Assume for $x$ tending to infinity}
                                      \[ \ds \sum_{n \leq x} b(n) = M(x) + \text{Er} (x), \]\textit{where}
                                     \begin{gather*}
                                        M(x) := \alpha x^2 \quad \textit{($\alpha$ is a certain complex number}),
                                        \\
                                        \text{Er} (x) := \sum_{n \leq x} b(n) - M(x).
                                     \end{gather*}
                                     \textit{Now, we consider the following Volterra integral equation of second type}
                                     \begin{equation*}
                                        F_1 (x) - \int_{0}^x F_1 (t)\frac{dt}{t} = \text{Er} (x) \quad (x \geq 0).
                                     \end{equation*}\textit{Then, for every complex number $A$, the function}
                                     \begin{equation*}
                                        F_1 (x) = (f_1 (x) + A)x  \quad (x \geq 0),
                                     \end{equation*}\textit{is a solution of the integral equation and these exhaust all solutions. Here,}
                                     \begin{equation*} 
                                        f_1 (x) = -\sum_{n=1}^\infty \frac{a(n)}{n}\left\{ \frac{x}{n} \right\}
                                     \end{equation*}\textit{for every $x \geq 0$.}

In the above previous research by {\K}, the arithmetical functions $a(n) \text{ and } b(n)$ in this conjecture are $\mu(n) \text{ and } \varphi(n)$ respectively, and all the hypothesis are satisfied. As for the error term $\text{Er}(x)$, we have a bound similar to $\text{Er}(x) = o(x^2)$ as $x$ tends to infinity in mind. Assuming certain appropriate hypothesis, in the present paper we prove this conjecture under a certain additional condition.

\textsc{\large Theorem.} \textit{Let $\{ a (n) \}$ be a complex-valued arithmetical function for which the series}
                                  \begin{equation*}
                                     \sum_{n=1}^\infty \frac{a (n)}{n^2}   \tag{1.5}
                                  \end{equation*}\textit{is convergent with the sum $2\alpha$, where $\alpha$ is an arbitrary complex number. Let $\{ b (n) \}$ be the arithmetical function defined by} 
                                  \begin{equation}
                                     b (n) = \sum_{d|n} a (d)\frac{n}{d}.   \tag{1.6}
                                  \end{equation}\textit{Assume for $x$ tending to infinity}
                                  \begin{equation*}
                                     \sum_{n \leq x} b (n) = M(x) + \text{Er}(x),   \tag{1.7} 
                                  \end{equation*}\textit{where}
                                  \begin{gather*}
                                     M(x) := \alpha x^2,   \tag{1.8}
                                     \\
                                     \text{Er}(x) := \sum_{n \leq x} b (n) -M(x).   \tag{1.9}
                                  \end{gather*}\textit{Now, we consider the following Volterra integral equation of second type}
                                   \begin{equation*}
                                      F_1 (x) - \int_{0}^x F_1 (t)\frac{dt}{t} = \text{Er}(x) \quad (x \geq 0).  \tag{1.10}
                                   \end{equation*}\textit{Then, for every complex number $A$, the function}
                                   \begin{equation*}
                                      F_1 (x) = (f_1 (x) + A)x  \quad (x \geq 0),   \tag{1.11} 
                                   \end{equation*}\textit{is a solution of the integral equation (1.10) and these exhaust all solutions of (1.10). Here, }
                                   \begin{equation*} 
                                      f_1 (x) = -\sum_{n=1}^\infty \frac{a(n)}{n}\left\{ \frac{x}{n} \right\}   \tag{1.12}
                                   \end{equation*}\textit{for every $x \geq 0$.}

As usual, if we say a function $F_1 $ is a solution of (1.10), then we implicitly assume that the integral in (1.10) exists in the sense that the limit
\begin{equation}
   \lim_{\epsilon \to 0+} \int_{\epsilon}^x F_1 (t)\frac{dt}{t}   \tag{1.13}
\end{equation}exists. We use the same convention throughout this paper. The formula (1.11) is a generalization of the result of Kaczorowski and Wiertelak~\cite{Kac and Wie}, except that the term $Ax$ is missing. Also,  the function $f_1 (x)$ is locally bounded. In fact, by the condition of theorem 
\begin{align*}
   f_1 (x)
   &= -\sum_{n=1}^\infty \frac{a (n)}{n}\left\{ \frac{x}{n} \right\}
   \\
   &= -\sum_{n=1}^\infty \frac{a (n)}{n} \left( \frac{x}{n} - \left[ \frac{x}{n} \right] \right)
   \\
   &= -x\sum_{n=1}^\infty \frac{a (n)}{n^2} + \left( \sum_{n \leq x}\frac{a (n)}{n}\left[ \frac{x}{n} \right] + \sum_{n > x}\frac{a (n)}{n} \left[ \frac{x}{n} \right] \right)
   \\
   &= -2\alpha x + \sum_{n \leq x}\frac{a (n)}{n}\left[ \frac{x}{n} \right].
\end{align*}

Here, we mention a previous research related to~\cite{Kac and Wie}. For $x\geq0$ let us put 
\begin{equation*}
   g(x) = \sum_{n=1}^\infty \mu(n)\left\{ \frac{x}{n} \right\}^{2}.  \tag{1.14}
\end{equation*}Then, $E(x)$ in (1.1) can be splitted as follows : 
\begin{equation*}
   E(x) = E^{\text{AR}}(x) + E^{\text{AN}}(x),   \tag{1.15}
\end{equation*}where
\begin{equation*}
   E^{\text{AR}}(x) = xf(x) \quad \text{and} \quad E^{\text{AN}}(x) = \frac{1}{2}g(x) + \frac{1}{2}   \tag{1.16}
\end{equation*}with $f(x)$ and $g(x)$ given by (1.3) and (1.14) respectively. ($E^{\text{AR}}(x) \text{ is the case } A \hspace{-0.1cm} = 0$ in (1.4).) We call $E^{\text{AR}}(x)$ and $E^{\text{AN}}(x)$ the \textit{arithmetic} and the \textit{analytic part} of $E(x)$ respectively. {\K} showed the $\Omega$-estimates for each of $E^{\text{AR}}(x)$ and $E^{\text{AN}}(x)$  (see~\cite{Kac and Wie},~\cite{Mon}). The decomposition (1.15) itself can obtain from certain calculations starting with  
\begin{equation*}
   \sum_{n \leq x} \varphi(n) = \frac{1}{2}\sum_{d \leq x} \mu(d)\left[ \frac{x}{d} \right]\left( \left[ \frac{x}{d} \right] + 1 \right).   \tag{1.17}
\end{equation*}

\begin{center}
\textbf{2.Related considerations for the twisted Euler $\varphi$-function}
\end{center}
\vspace{0.1cm}

{\K} obtained a better decomposition for the remainder term in the asymptotic formula for a generalization of the Euler totient function ({see~\cite{kac and wie}) : For a non-principal real Dirichlet character $\chi \hspace{0.1cm} (\text{mod }q), q>2,\text{ let } \varphi(n,\chi)$ denote the twisted Euler $\varphi$-function 
\vspace{-0.1cm}
\[ \ds \varphi(n,\chi) = n\prod_{p|n} \left( 1 - \frac{\chi(p)}{p} \right). \]

{\K} made a similar consideration to~\cite{Kac and Wie} for the remainder term in the asymptotic formula of the above twisted Euler $\varphi$-function.

Let
\vspace{-0.1cm}
\begin{equation*}
   E(x,\chi) = \sum_{n \leq x} \varphi(n,\chi) - \frac{x^2}{2L(2,\chi)}   \tag{2.1}
\end{equation*}
\vspace{-0.1cm}
and
\vspace{-0.1cm}
\begin{equation*}
   E_1 (x,\chi) = \begin{cases}
                         E(x,\chi) \quad (x \notin \mathbb{N}),
                         \\
                         \frac{1}{2}(E(x-0,\chi)+E(x+0,\chi)) \quad (\text{otherwise})   \tag{2.2}
                      \end{cases}
\end{equation*}be the corresponding error term. Here, as usual, $L(s,\chi)$ denotes the Dirichlet $L$-function associated to $\chi$. It is easy to see that $E(x,\chi) = O(x\log x)$ for $x \geq 2$.

Let $s(x)$ be the saw-tooth function
\vspace{-0.1cm}
\begin{equation*}
   s(x) = \begin{cases}
               0 \quad (x \in \mathbb{Z})            
               \\
               \frac{1}{2} - \{ x \} \quad (\text{otherwise}).   \tag{2.3}
            \end{cases}
\end{equation*}

We write for $x \geq 0$
\vspace{-0.1cm}
\begin{gather*}
   f(x,\chi) = \sum_{d=1}^\infty \frac{\mu(d)\chi(d)}{d}s\left( \frac{x}{d} \right),   \quad \tag{2.4}
   \\
   g(x,\chi) = \sum_{d=1}^\infty \mu(d)\chi(d) \left\{ \frac{x}{d} \right\}\left( \left\{ \frac{x}{d} \right\} - 1 \right).  \quad  \tag{2.5}
\end{gather*}
\vspace{-0.1cm}
Then the solution of the following Volterra integral equation of second type
\[ F(x,\chi) - \int_{0}^\infty K(x,t)F(t,\chi)dt = E_1 (x,\chi) \quad (x \geq 0), \]where 
\vspace{-0.1cm}
\begin{equation*}
   K(x,t) = \begin{cases}
                 1/t \quad (0 < t \leq x),
                 \\
                  0  \quad (0 \leq x < t), 
              \end{cases}
\end{equation*}is the function
\vspace{-0.1cm}
\begin{equation*}
   F(x,\chi) = (f(x,\chi)+A)x,   \tag{2.6}
\end{equation*}where $A$ is an arbitrary constant. (In the paper~\cite{kac and wie}, the unique solution is $F(x,\chi) = xf(x,\chi)$, but the comments just after (1.4) should also be applied here.)

Moreover, when $A = 0$ in (2.6), for $x \geq 0$ 
\vspace{-0.1cm}
\begin{equation*}
   E_1 (x,\chi) = E^{\text{AR}}(x,\chi) + E^{\text{AN}}(x,\chi),  \tag{2.7}
\end{equation*}where
\vspace{-0.1cm}
\begin{equation*}
   E^{\text{AR}}(x,\chi) = xf(x,\chi) \quad \text{and} \quad E^{\text{AN}}(x,\chi) = \frac{1}{2}g(x,\chi)   \tag{2.8}
\end{equation*}with $f(x,\chi)$ and $g(x,\chi)$ given by (2.4) and (2.5) respectively.  

\textsc{Remark}
(1) The function (2.4) is a better solution of the Volterra integral equation (1.10) than (1.12), namely, the decomposition (2.7) coincides with the decomposition (1.15). For simplicity, let us forget about  the error term (2.2) and consider $E(x,\chi)$. Let
\begin{equation*}
   s(x) := \frac{1}{2} - \{ x \},   \tag{2.9}
\end{equation*}which is not the same (2.3), but (2.3) is just the normalized version of (2.9). Then we have (after removing normalization)
\begin{equation*}
   f(x,\chi) := \sum_{d=1}^\infty \frac{\mu(d)\chi(d)}{d}s\left( \frac{x}{d} \right) 
                 = -\sum_{d=1}^\infty  \frac{\mu(d)\chi(d)}{d}\left\{ \frac{x}{d} \right\} + \frac{1}{2}\sum_{d=1}^\infty \frac{\mu(d)\chi(d)}{d}  \tag{2.10}
\end{equation*}and
\begin{align*}
   g(x,\chi)
   &= \sum_{d=1}^{\infty} \mu(d)\chi(d) \left\{ \frac{x}{d} \right\} \left(  \left\{ \frac{x}{d} \right\} - 1 \right)
   \\
   &= \sum_{d=1}^{\infty} \mu(d)\chi(d) \left\{ \frac{x}{d} \right\}^2 
     - \sum_{d=1}^{\infty}  \mu(d)\chi(d) \left\{ \frac{x}{d} \right\}.
\end{align*} If $\chi$ is trivial, we have
\begin{equation*}
   \frac{1}{2}\sum_{d=1}^\infty \frac{\mu(d)\chi(d)}{d} = \frac{1}{2}\sum_{d=1}^\infty \frac{\mu(d)}{d} = 0 \end{equation*}and
\begin{align*}
   - \sum_{d=1}^\infty  \mu(d)\chi(d) \left\{ \frac{x}{d} \right\}
   &=  - \sum_{d=1}^\infty  \mu(d) \left\{ \frac{x}{d} \right\}
   \\
   &= \sum_{d=1}^\infty \mu(d)\left[ \frac{x}{d} \right] - \sum_{d=1}^\infty  \mu(d) \cdot \frac{x}{d}
    \\
    &= \sum_{d \leq x} \mu(d)\left[ \frac{x}{d} \right]
    \\
    &= 1.
\end{align*}Thus, for the trivial character, we have
\[  f(x,\chi) = f(x) \quad  \text{and} \quad g(x,\chi) = g(x) + 1.  \]However, if we use the main theorem to $a (n)=\mu(n)\chi(n)$, then the arithmetic part becomes
\[ \displaystyle -x\sum_{d=1}^\infty \frac{\mu(d)\chi(d)}{d}\left\{ \frac{x}{d} \right\} \]as in (1.12) and so this differs form (2.10). This discrepancy happened, indeed, because of the missing term $Ax$. The remaining term
\[ \displaystyle  x \cdot \frac{1}{2}\sum_{d=1}^\infty \frac{\mu(d)\chi(d)}{d} \]of $x \cdot f(x,\chi)$ is just of the form $Ax$ and so such term is allowed for the general solution of the integral equation (1.10). As a summary, the function (2.4) is a better representative of the solutions of the Volterra integral equation (1.10). However, this representative is that the solution (2.4) is well-defined only when the series  
\begin{equation*}
   \sum_{d=1}^\infty \frac{\mu(d)\chi(d)}{d}  \tag{2.11}  
\end{equation*}converges.

(2) In the paper ~\cite{BP}, the asymptotic formulas related to the divisor function $\sigma(n)$ and the Euler totient function $\varphi(n)$ are mentioned (see~\cite{BP}, p61,62). By introducing the weight $n$ into the formula by using partial summation and using Theorem 4 of ~\cite{BP}, for example, we get  
      \[ \displaystyle \sum_{n \leq x} n^{1-\alpha}\sigma^\alpha (n) = Cx^2 + \sum_{r=0}^{[\alpha]-1}C_r x(\log x)^{\alpha-r-1} + O(x(\log x)^\frac{2|\alpha|}{3} (\log\log x)^\frac{4|\alpha|}{3}). \]
     $b (n)=n^{1-\alpha}\sigma^\alpha (n)$ fits into the setting of the main theorem. Thus, we have second main terms
     \[  \displaystyle \sum_{r=0}^{[\alpha]-1}C_r x(\log x)^{\alpha-r-1} \]  provided $\alpha$ is sufficiently large. These second main terms cannot be seen in the result of ~\cite{Kac and Wie}.

\begin{center}
\textbf{3.Proof of Theorem}
\end{center}
Now we prove the theorem. We define the auxiliary function for $x \geq 0$ by
\vspace{-0.1cm}
\begin{equation*}
   R(x) = \text{Er}(x) - xf_1  (x).   \tag{3.1}
\end{equation*}First, we prepare the following two lemmas.

\textsc{LEMMA1.} \textit{For all positive}   $x$,
                         \begin{equation*}
                            R(x) = -\int_{0}^x f_1 (t)dt.  \tag{3.2}
                         \end{equation*}
\textit{Proof.} Let us observe that $R(x)$ is a continuous function. For $x = 0$ and for positive $x$ which is not an integer, it is evident. Let $N$ be a positive integer. By splitting the series (1.12) at $N$, and considering the limit $\left\{ (N+x)/n \right\}$ as $x$ tending to $0$, we see that
\begin{align*}
   f_1 (N+0)
   &= -\sum_{n=1}^\infty \frac{a (n)}{n}\left\{ \frac{N+0}{n} \right\},
   \\
   f_1 (N-0)
   &= -\sum_{n=1}^\infty \frac{a (n)}{n}\left\{ \frac{N-0}{n} \right\}.
\end{align*}Since
\begin{equation*}
   \left\{ \frac{N+0}{n} \right\} - \left\{ \frac{N-0}{n} \right\}
   =\begin{cases}
        0  \quad \text{($n \nmid N$)}
        \\
        -1 \quad \text{($n \mid N$)}
    \end{cases}
\end{equation*}(see~\cite{Kac and Wie}, P2691), we have
\vspace{-0.1cm}
\[ \ds f_1 (N+0) - f_1 (N-0) = \sum_{n|N} \frac{a (n)}{n} = \frac{b (N)}{N}. \]
Therefore
\begin{eqnarray*}
   R(N+0) - R(N-0)
\hspace{-0.25cm}   &=& \hspace{-0.25cm} (\text{Er}(N+0) - \text{Er}(N-0)) - N( f_1 (N+0) - f_1 (N-0) )
   \\
   &=& \hspace{-0.25cm} b (N) - N \cdot \frac{b (N)}{N}
   \\
   &=& \hspace{-0.25cm} 0,
\end{eqnarray*}and hence $R(N-0) = R(N+0) =R(N)$. 

Let $x$ be positive and not an integer. Take derivatives of the both sides of (3.1). Since $x$ is not a positive integer, we have $\text{Er}^\prime (x) = -M^\prime (x) = -2\alpha x$. Therefore we have
\[ R^\prime (x) = -2\alpha x - f_1 (x)-xf_{1}^\prime (x). \]For $x$ which is positive and not an integer, we have $\{ x/n \}^\prime = 1/n$ (see~\cite{Kac and Wie}, p2691). Considering the hypothesis on the series (1.5), differentiating term by term we obtain
\vspace{-0.1cm}
\[ \ds f_{1}^\prime (x) = -\sum_{n=1}^\infty \frac{a (n)}{n} \cdot \frac{1}{n} = -2\alpha. \]Consequently, we have
\vspace{-0.1cm}
\[ R^\prime (x) = -f_{1}(x) \]for $x$ which is positive and not an integer. Because of $R(0) = 0$ and the continuity of $R(x)$, we have (3.2) for all positive $x$. $\square$

\textsc{LEMMA2.} \textit{Let $G$ be a complex-valued function defined on $[0, \infty)$ satisfying}
                         \begin{equation}
                            \int_{0}^x |G(t)|\frac{dt}{t} < +\infty   \tag{3.3}
                         \end{equation}\textit{and the integral equation}
                         \begin{equation}
                            G(x) - \int_{0}^x G(t)\frac{dt}{t} = 0   \tag{3.4}
                         \end{equation}\textit{for all $x \geq 0$. Then we have}
                         \begin{equation}
                            G(x) = Ax    \tag{3.5}
                        \end{equation}\textit{for some complex number $A$.}

\textit{Proof.} It is obvious that (3.5) satisfies (3.3) and (3.4) for all $x \geq 0$. Conversely, take a function $G(x)$ arbitrarily satisfying (3.3) and (3.4) for all $x \geq 0$. By (3.3) and (3.4), we see that 
\[ G(x) = \int_{0}^x G(t)\frac{dt}{t} \] is a continuous function on $[0,+\infty)$. Thus, using integral equation again and using the fundamental theorem of calculus, we see that $G(x)$ is continuously differentiable on $(0,+\infty)$. By taking the derivative of (3.4), we have
\[ G^\prime (x) = \frac{G(x)}{x} \quad (x>0). \]Thus, we have $G(x) = Ax$ for $x > 0$ for some $A$ and by the continuity this holds for $x \geq 0$. $\square$

\textsc{Proof of Theorem.} Let a function $F_1 (x)$ be the solution of the Volterra integral equation of second type (1.10) satisfying the condition (1.13). Using (3.1) and (1.10), from (3.2) we have 
\begin{equation*} 
   \int_{0}^x t^{-1}( F_1 (t) - tf_1 (t) )dt = F_1 (x) - xf_1 (x) \quad (x \geq 0).   \tag{3.6}
\end{equation*}Now we put 
\vspace{-0.1cm}
\begin{equation*}
   G(x) := F_1 (x) - xf_1 (x).   \tag{3.7} 
\end{equation*}Then, the equation (3.6) yields
\begin{equation*} 
   \int_{0}^x t^{-1}G(t)dt = G(x) \quad (x \geq 0).    \tag{3.8}
\end{equation*}Using Lemma 2, we must have (3.5). By substituting into (3.7), we have the solution (1.11).

Conversely, if we assume that $F_1 (x)$ is a function of type (1.11). Then,
\begin{align*}
   F_1 (x) - \int_{0}^x F_1 (t)\frac{dt}{t}
   &= (f_1 (x) + A)x - \int_{0}^x (f_1 (t) + A)dt
   \\
   &= (f_1 (x) + A)x - \int_{0}^x f_1 (t)dt - Ax
   \\
   &= xf_1 (x) - \int_{0}^x f_1 (t)dt.
\end{align*}Using (3.1) and (3.2), 
\begin{align*}
   xf_1 (x) - \int_{0}^x f_1 (t)dt
   &= xf_1 (x) + R(x)
   \\ 
   &= xf_1 (x) +  \text{Er}(x) - xf_1 (x)
   \\
   &= \text{Er}(x).
\end{align*}Therefore, the function $F_1 (x)$ of type (1.11) is the solution of the integral equation (1.10) for all $x \geq 0$. Since the function $f_1 (x)$ is a locally bounded as noted in section 1, and $A$ is a constant, it is clear that the function $F_1 (x)$ satisfies the condition (1.13). The completes the proof. $\square$

\textsc{Remark.} It is possible to deduce our theorem from a known result. In fact, the entry 2.1.50 of ~\cite{P and M} with
\[ A := -1, \> \lambda := 0, \> \mu := -1, \> f(x) := \text{Er} (x), \] which satisfies $\lambda + \mu + 1 = 0$, gives the solution (1.11) of the Volterra integral equation of second type (1.10). We can show the following fact : 

\textit{Let $E_1$ be a complex-valued function defined on $[0, +\infty)$ satisfying}
\begin{equation*}
   \int_{0}^x |E_1 (t)|\frac{dt}{t^2} < +\infty   \tag{3.9}
\end{equation*}\textit{for all $x \geq 0$. Then, the function}
\begin{equation}
   F_2 (x) := E_1 (x) + x\int_{0}^x E_1 (x)\frac{dt}{t^2}   \tag{3.10}
\end{equation}\textit{well-defined on $[0,+\infty)$ satisfies the Volterra integral equation of second type}
\begin{equation}
   F_2 (x) - \int_{0}^x F_2 (t)\frac{dt}{t} = E_1 (x)   \tag{3.11}
\end{equation}\textit{and satisfies}
\begin{equation*}
   \int_{0}^x |F_2 (t)|\frac{dt}{t} < +\infty   \tag{3.12}
\end{equation*}\textit{for all $x \geq 0$.}

It follows immediately that the special solution of the integral equation (1.10) is expressed using (1.12) by using (3.10). The general solution of (3.11) can be expressed as follows : 

\textit{Let $E_1$ be a complex-valued function defined on $[0, +\infty)$ satisfying the condition (3.9) for all $x \geq 0$. Then, the complex-valued functions $F_3 $ defined on $[0, +\infty)$ satisfying the condition (3.12) and the Volterra integral equation of second type (3.11) are given by}
\begin{equation*}
   F_3 (x) = E_1 (x) +x\int_{0}^x E_1 (t)\frac{dt}{t^2} + Ax    \tag{3.13}
\end{equation*}\textit{with complex number $A$.}

However, in this paper, we choose a self-contained method how to construct the solution concretely not using the formula in the entry 2.1.50 of ~\cite{P and M} . 

\vspace{0.1cm}

\textbf{Acknowledgments.}   The author expresses his sincere gratitude to Prof. \\ Kohji Matsumoto, Dr. Yuta Suzuki and Dr. Wataru Takeda for giving me various advice in writing this paper, Prof. Isao Kiuchi for introducing the reference~\cite{Kac and Wie}, and Dr. Shota Inoue for providing the conjecture on page 2. Also, the author appreciates valuable comments from the anonymous referee. In particular, the author referred to comments from the referee about the remark of section 2.

\bigskip
　　　　　　　　　　　　　　Hideto Iwata 
                                                  
　　　　　　　　　　　　　　Graduate School of Mathematics

　　　　　　　　　　　　　　 Nagoya University
                                                   
　　　　　　　　　　　　　　 Furocho, Chikusa-ku,
                                                   
　　　　　　　　　　　　　　 Nagoya, 464-8602, Japan.

　　　　　　　　　　　　　　 \small{e-mail:d18001q@math.nagoya-u.ac.jp}

\begin{thebibliography}{99}
\bibitem{BP} 
U. Balakrishnan, Y.-F.  S. P\'{e}termann,
\newblock The Dirichlet series of $\zeta(s)\zeta^\alpha (s+1)f(s+1)$ : On an error term associated with its coefficients,
\newblock Acta Arithmetica \textbf{75} (1996), 39-69. 

\bibitem{Kac and Wie} 
J. Kaczorowski, K. Wiertelak, 
\newblock Oscillations of the remainder term related to the Euler totient function, 
\newblock J. Number Theory \textbf{130} (2010) 2683-2700. 

\bibitem{kac and wie}
J. Kaczorowski, K. Wiertelak, 
On the sum of the twisted Euler function,  
\newblock Int. J. Number Theory \textbf{8} (7) (2012) 1741-1761. 

\bibitem{Mon}
H. L. Montgomery, 
Fluctuations in the mean of Euler's phi function,  
\newblock Proc. Indian Acad. Sci. Math. Sci. \textbf{97} (1-3) (1987) 239-245.

\bibitem{P and M}
A. D. Polyanin and A. V. Manzhirov,
Handbook of Integral Equations 2nd edition, Chapman \& Hall / CRC, (2008)
\end{thebibliography}
\end{document}